\documentclass[twoside,english,final, 1p, times]{elsarticle}
\usepackage[T1]{fontenc}
\usepackage[latin9]{inputenc}
\usepackage{geometry}
\usepackage{fancyhdr}
\usepackage{color}
\usepackage{babel}
\usepackage{amsthm}
\usepackage{amsmath}
\usepackage{amssymb}
\usepackage{graphicx}
\usepackage{esint}
\usepackage[unicode=true,
 bookmarks=true,bookmarksnumbered=false,bookmarksopen=false,
 breaklinks=true,pdfborder={0 0 1},backref=false,colorlinks=true]
 {hyperref}
\hypersetup{pdftitle={Riemann Zeta},
 pdfauthor={Alois Pichler}}

\makeatletter

\numberwithin{figure}{section}
\theoremstyle{plain}
\newtheorem{thm}{\protect\theoremname}
  \theoremstyle{remark}
  \newtheorem{rem}[thm]{\protect\remarkname}

\journal{arXiv.org}
\renewcommand{\Re}{\operatorname{Re}}

\makeatother

  \providecommand{\remarkname}{Remark}
\providecommand{\theoremname}{Theorem}

\begin{document}

\title{Strategies To Evaluate The Riemann Zeta Function}

\author{Alois Pichler}

\address{Department of Statistics and Operations Research, University of Vienna,
Austria, 1010 Wien, Universitätsstraße~5}

\ead{alois.pichler@univie.ac.at}
\begin{abstract}
This paper continues a series of investigations on converging representations
for the Riemann Zeta function. We generalize some identities which
involve Riemann's zeta function, and moreover we give new series and
integrals for the zeta function. The results originate from attempts
to extend the zeta function by classical means on the complex plane.
This is particularly of interest for representations which converge
rapidly in a given area of the complex plane, or for the purpose to
make error bounds available. \end{abstract}
\maketitle

\section{Introduction and Definitions}

\subsection{Extensions To The Complex Plane}

The Riemann Zeta function is classically defined as 
\begin{equation}
\zeta\left(s\right)=\sum_{i=1}^{\infty}\frac{1}{i^{s}},\label{eq:Zeta}
\end{equation}
which is a convergent representation for $\Re s>1$. An immediate
extension to $\Re s>0$ is given by the variant 
\[
\zeta\left(s\right)=\frac{1}{1-2^{1-s}}\sum_{i=1}^{\infty}\frac{\left(-1\right)^{i-1}}{i^{s}},
\]
which results from separating summands of the form $\frac{1}{\left(2i\right)^{s}}$
and $\frac{1}{\left(2i+1\right)^{s}}$. This alternating series converges
for $\Re s>0$ and reveals that $\zeta$ has a pole with residue 1
at $s=1$, as $\sum_{i=1}\frac{\left(-1\right)^{i-1}}{i}=\ln2$.

Probably less well-known is the similar representation 
\begin{eqnarray*}
\zeta\left(s\right) & = & \frac{1}{1-3^{1-s}}\sum_{i=1}^{\infty}\left(\frac{1}{\left(3i-2\right)^{s}}+\frac{1}{\left(3i-1\right)^{s}}-\frac{2}{\left(3i\right)^{s}}\right)\\
 & = & \frac{1}{1-3^{1-s}}\left(1+\sum_{i=1}^{\infty}\left(\frac{1}{\left(3i-1\right)^{s}}-\frac{2}{\left(3i\right)^{s}}+\frac{1}{\left(3i+1\right)^{s}}\right)\right),
\end{eqnarray*}
where the latter converges as well for $\Re s>-1$.

Other interesting representations, as they seem to not appear in the
literature and which extend the zeta function even further, include
\begin{itemize}
\item $\left(1-\frac{1}{2^{s}}-\frac{2}{4^{s}}\right)\zeta\left(s\right)=1+\sum_{i=1}\left(\frac{1}{\left(4i-1\right)^{s}}-\frac{2}{\left(4i\right)^{s}}+\frac{1}{\left(4i+1\right)^{s}}\right)$
for $\Re s>-1$,
\item $\left(1-\frac{5}{2^{s}}+\frac{5}{3^{s}}-\frac{1}{6^{s}}\right)\zeta\left(s\right)=\sum_{i=0}\left(\frac{1}{\left(6i+1\right)^{s}}-\frac{4}{\left(6i+2\right)^{s}}+\frac{6}{\left(6i+3\right)^{s}}-\frac{4}{\left(6i+4\right)^{s}}+\frac{1}{\left(6i+5\right)^{s}}\right)$
and
\item $\left(4-\frac{5}{2^{s}}-\frac{4}{3^{s}}-\frac{1}{6^{s}}\right)\zeta\left(s\right)=4-\frac{1}{2^{s}}-\sum_{i=1}\left(\frac{1}{\left(6i-2\right)^{s}}-\frac{4}{\left(6i-1\right)^{s}}+\frac{6}{\left(6i\right)^{s}}-\frac{4}{\left(6i+1\right)^{s}}+\frac{1}{\left(6i+2\right)^{s}}\right)$
for $\Re s>-3$;
\end{itemize}
they are all verified in a similar way as the initial formula.

\subsection{Euler-Maclaurin}

To further extend the Zeta function to complex arguments with even
smaller real part one may apply Euler-Maclaurin's summation formula,
which states that 
\[
\sum_{i=1}^{n-1}f\left(i\right)=\int_{1}^{n}f\left(x\right)\mathrm{d}x+\sum_{j=1}^{\Delta}\frac{B_{j}}{j!}\left(f^{\left(j\right)}\left(n\right)-f^{\left(j\right)}\left(1\right)\right)-\left(-1\right)^{\Delta}\int_{1}^{n}f^{\left(\Delta\right)}\left(x\right)\frac{B_{\Delta}\left(\left\{ x\right\} \right)}{\Delta!}\mathrm{d}x;
\]
$f$ is a function with (continuous) derivatives $f^{\left(i\right)}$
up to order $\Delta$, $\left\{ x\right\} =x-\left\lfloor x\right\rfloor $
is the fractional part of $x$, $B_{i}$ the Bernoulli number and
$B_{i}\left(x\right)=\sum_{j=0}^{i}{i \choose j}B_{i-j}x^{j}$ the
Bernoulli polynomial (cf. \citep{abramowitz+stegun}).

The statement may be applied to \eqref{eq:Zeta}, giving thus after
some rearrangements the identity 
\begin{equation}
\zeta\left(s\right)=\sum_{i=1}^{n-1}\frac{1}{i^{s}}+\sum_{j=0}^{\Delta}\frac{B_{j}}{n^{s+j-1}}\frac{{1-s \choose j}}{s-1}-\left(-1\right)^{\Delta}{-s \choose \Delta}\cdot\int_{n}^{\infty}\frac{B_{\Delta}\left(\left\{ x\right\} \right)}{x^{s+\Delta}}\mathrm{d}x\label{eq:Euler-Maclaurin}
\end{equation}

(${z \choose n}=\frac{\Gamma\left(z+1\right)}{n!\Gamma\left(z-n+1\right)}=\frac{z}{n}\cdot\frac{z-1}{n-1}\cdots\frac{z-n+1}{1}$
is the (generalized) binomial coefficient.)

The error term $\int_{n}^{\infty}\frac{B_{\Delta}\left(\left\{ x\right\} \right)}{x^{s+\Delta}}\mathrm{d}x$
often is inconvenient, as it involves an integral of a function, which
is not smooth enough. We shall elaborate different expressions below,
replacing this error term typically by other series or by different
integrals; as an example recall the integral 
\begin{equation}
\zeta\left(s\right)=1+\frac{1}{s-1}\sum_{k=0}^{\Delta}{k+s-2 \choose k}B_{k}+\frac{1}{\Gamma\left(s\right)}\int_{0}^{\infty}\frac{1}{e^{t}-1}-\sum_{k=0}^{\Delta}\frac{B_{k}}{k!}t^{k-1}e^{-t}t^{s-1}\mathrm{d}t\label{eq:Exp-Integral}
\end{equation}
(cf. \citep{Choi2010} or \citep{Srivastava2001}).

This construction plays even a prominent role in the famous algorithm
of Odlyzko and Schönhage's to evaluate $\zeta$ on the critical strip,
so some remarks seem to be appropriate (cf. \citep{Odlyzko+Schoenhage},
\citep{Cohen1992} and \citep{Petermann2007}).
\begin{rem}
Identity \eqref{eq:Euler-Maclaurin} is frequently given in the form
\[
\zeta\left(s\right)=\sum_{i=1}^{n}\frac{1}{i^{s}}-\sum_{j=0}^{\Delta}\frac{B_{j}}{n^{s+j-1}}\frac{{s+j-2 \choose j}}{1-s}-\left(-1\right)^{\Delta}{-s \choose \Delta}\cdot\int_{n}^{\infty}\frac{B_{\Delta}\left(\left\{ x\right\} \right)}{x^{s+\Delta}}\mathrm{d}x,
\]
which is equivalent because $B_{2j+1}=0$ (except for $B_{1}=-\frac{1}{2}$)
and ${1-s \choose j}=\left(-1\right)^{j}{s+j-2 \choose j}$.
\end{rem}

\begin{rem}
As $B_{\Delta}$ is a polynomial, $B_{\Delta}\left(\left\{ x\right\} \right)$
is uniformly bounded for $x\in\mathbb{R}$, and it even holds that
$\left|B_{\Delta}\left(x\right)\right|\le\frac{2\,\Delta!}{\left(1-2^{1-\Delta}\right)\left(2\pi\right)^{\Delta}}$
for $x\in\left[0,1\right]$ (cf. \citep{Lehmer1940}). This allows
to conclude that%
\footnote{Note again that $B_{2j+1}=0$, so the order of convergence in the
next statement actually is $\mathcal{O}\left(\frac{1}{n^{s+\Delta}}\right)$,
whenever $\Delta$ is even.%
} 
\[
\zeta\left(s\right)=\sum_{i=1}^{n-1}\frac{1}{i^{s}}-\sum_{j=0}^{\Delta}\frac{B_{j}}{n^{s+j-1}}\frac{{1-s \choose j}}{1-s}+\mathcal{O}\left(\frac{1}{n^{s+\Delta-1}}\right),
\]
and whence, $\zeta$ is analytic in the entire complex plane except
for $s=1$, where is a pole with residue $1$. The argument, however,
can be refined for even more explicit error bounds and rates (cf.
\citep{Cohen1992}). But as a general patter the expression converge,
provided that $\Re s>1-\Delta$.
\end{rem}

\begin{rem}
Notice as well that ${-s \choose \Delta}=0$ for $s$ a negative integer
and $\Delta$ big enough, which enables to recover 
\[
\zeta\left(-n\right)=-\frac{B_{n+1}}{n+1}
\]
for $n\in\left\{ 1,2,\dots\right\} $.
\end{rem}

\begin{rem}
The partial series is frequently denoted 
\[
H_{n}\left(s\right):=\sum_{i=1}^{n}\frac{1}{i^{s}},
\]
which is the generalized harmonic number of order $n$ of $s$.
\end{rem}

\begin{rem}
For future reference recall that 
\begin{equation}
H_{n}\left(-k\right)=\sum_{i=1}^{n}i^{k}=\frac{B_{k+1}\left(n+1\right)-B_{k+1}}{k+1}\label{eq:Faulhaber}
\end{equation}
in the notation just introduced for natural numbers $k\in\left\{ 0,1,2,\dots\right\} $
-- this is sometimes referred to as \emph{Faulhaber's formula}. 
\end{rem}

\section{Results And Extension To Various Directions}

It turns out that a lot of different relations for the Riemann Zeta
function again involve the expression $\sum_{j=0}^{\Delta}B_{j}\frac{{1-s \choose j}}{s-1}$.
This is natural, as we have seen that this expression describes the
asymptotic behavior of $\sum_{i=1}^{n}i^{-s}$.

So does the following extension of an identity which is sometimes
referred to as \emph{Stark's formula}; an extension to Hurwitz Zeta
function is known under the name \emph{Stark-Keiper formula}:
\begin{thm}[Generalization of Stark's formula]
For $s\neq1$ and any $\Delta\in\left\{ 1,2,\dots\right\} $, 
\begin{equation}
\zeta\left(s\right)=1+\frac{1}{s-1}\sum_{k=0}^{\Delta-1}{s+k-2 \choose k}B_{k}-\frac{1}{s-1}\sum_{i=\Delta}\left(\zeta\left(s+i\right)-1\right){s+i-1 \choose i+1}\sum_{k=0}^{\Delta-1}{i+1 \choose k}B_{k}.\label{eq:Stark}
\end{equation}
\end{thm}
\begin{rem}
The latter series converges for all $s\neq1$. As $\zeta\left(s+i\right)$
($i\ge\Delta$) are easily available even by means of formula \eqref{eq:Zeta}
provided that $\Re s>1-\Delta$, the statement makes $\zeta\left(s\right)$
available by different choices of $\Delta$ for $\Re s\le1$.

The particular statement for $\Delta=1$ (sometimes called \emph{Stark's
formula}, although it seems it was published first by Landau) is (for
example) contained in \citep{Titchmarsh}.\end{rem}
\begin{proof}
Notice first that 
\begin{alignat}{1}
1-\frac{1}{s-1}\sum_{i=1} & \left(\zeta\left(s+i\right)-1\right){s+i-1 \choose i+1}=\nonumber \\
= & 1-\frac{1}{s-1}\sum_{n=2}\sum_{i=1}\frac{1}{n^{s+i}}{s+i-1 \choose i+1}\label{eq:Stark-1}\\
= & =1-\frac{1}{s-1}\sum_{n=2}\sum_{i=2}\frac{\left(-1\right)^{i}}{n^{s+i-1}}{1-s \choose i}=1-\frac{1}{s-1}\sum_{n=2}\frac{1}{n^{s-1}}\left(\left(1-\frac{1}{n}\right)^{1-s}-1+\frac{1-s}{n}\right)\nonumber \\
= & 1-\frac{1}{s-1}\sum_{n=2}\frac{1}{\left(n-1\right)^{s-1}}-\frac{1}{n^{s-1}}+\frac{1-s}{n^{s}}=1-\frac{1}{s-1}\left(1+\left(1-s\right)\left(\zeta\left(s\right)-1\right)\right)\nonumber \\
= & \zeta\left(s\right)-\frac{1}{s-1},\nonumber 
\end{alignat}
which is the statement for $\Delta=1$.

The general result follows by induction on $\Delta$: To this end
assume that \eqref{eq:Stark} holds true and replace $s$ in \eqref{eq:Stark-1}
by $s+\Delta$. Subtracting appropriate multiples of \eqref{eq:Stark-1}
from \eqref{eq:Stark} to eliminate $\zeta\left(s+\Delta\right)$
reveals the assertion for $\Delta+1$, where additionally the identity
\[
\sum_{k=0}^{\Delta}{\Delta+1 \choose k}B_{k}=0
\]
($\Delta\ge1$) is involved. \end{proof}
\begin{thm}
For $\Delta\in\left\{ 1,2,\dots\right\} $ and $s\neq1$ the relation
\begin{equation}
\zeta\left(s\right)=\sum_{k=\Delta}{s+k-1 \choose k}\frac{\zeta\left(s+k\right)}{2^{k}}\sum_{j=0}^{\Delta-1}\frac{B_{j}}{2^{s}-2^{1-j}}\sum_{i=j}^{\Delta-1}{k \choose i}{i+1 \choose j}\frac{1}{i+1}\label{eq:Stark2-1}
\end{equation}
holds true.\end{thm}
\begin{rem}
The latter series converges for all $s\neq1$. As above, $\zeta\left(s+i\right)$
is easily available by means of \eqref{eq:Zeta} for $\Delta>1-\Re s$,
the statement whence again makes $\zeta\left(s\right)$ accessible
for arbitrary $\Re s\le1$.\end{rem}
\begin{proof}
The Statement for $\Delta=1$ is 
\begin{equation}
\zeta\left(s\right)=\frac{1}{2^{s}-2}\sum_{k=1}{s+k-1 \choose k}\frac{\zeta\left(s+k\right)}{2^{k}}\label{eq:Stark2}
\end{equation}
and sometimes called Ramaswami's formula. This is verified straight
forward, as 
\begin{alignat*}{1}
\sum_{k=1}{s+k-1 \choose k}\frac{\zeta\left(s+k\right)}{2^{k}} & =\sum_{i=1}\frac{1}{i^{s}}\sum_{k=1}{s+k-1 \choose k}\left(\frac{1}{2i}\right)^{k}\\
 & =\sum_{i=1}\frac{1}{i^{s}}\left(\left(1-\frac{1}{2i}\right)^{-s}-1\right)=\sum_{i=1}\frac{1}{i^{s}}\left(\left(\frac{2i}{2i-1}\right)^{s}-1\right)\\
 & =2^{s}\sum_{i=1}\left(\frac{1}{2i-1}\right)^{s}-\sum_{i=1}\frac{1}{i^{s}}=2^{s}\left(1-2^{-s}\right)\zeta\left(s\right)-\zeta\left(s\right)\\
 & =\left(2^{s}-2\right)\zeta\left(s\right).
\end{alignat*}
As for the general statement one may replace $s$ in \eqref{eq:Stark2}
by $s+\Delta$ and substitute this for the first summand in \eqref{eq:Stark2-1}.
This burns down then to the identity 
\[
\sum_{j=0}^{\Delta-1}\frac{B_{j}}{2^{s}-2^{1-j}}\sum_{i=j}^{\Delta-1}\frac{{\Delta \choose i}}{2^{s+\Delta}-2}\cdot\frac{{i+1 \choose j}}{i+1}=\sum_{j=0}^{\Delta}\frac{B_{j}}{2^{s}-2^{1-j}}\frac{{\Delta+1 \choose j}}{\Delta+1},
\]
which holds true because 
\begin{alignat*}{1}
\sum_{j=0}^{\Delta-1}\frac{B_{j}}{2^{s}-2^{1-j}} & \frac{1}{2^{s+\Delta}-2}\sum_{i=j}^{\Delta-1}{\Delta \choose i}{i+1 \choose j}\frac{1}{i+1}=\\
 & =\sum_{j=0}^{\Delta-1}\frac{B_{j}}{2^{s}-2^{1-j}}\frac{{\Delta+1 \choose j}}{\Delta+1}\frac{2^{1-j}-2^{1-\Delta}}{2^{s}-2^{1-\Delta}}\\
 & =\sum_{j=0}^{\Delta-1}B_{j}\frac{{\Delta+1 \choose j}}{\Delta+1}\left(\frac{1}{2^{s}-2^{1-j}}-\frac{1}{2^{s}-2^{1-\Delta}}\right)\\
 & =\sum_{j=0}^{\Delta}\frac{B_{j}}{2^{s}-2^{1-j}}\frac{{\Delta+1 \choose j}}{\Delta+1}
\end{alignat*}
and $\sum_{j=0}^{\Delta-1}B_{j}{\Delta+1 \choose j}=-\left(\Delta+1\right)B_{\Delta}$.\end{proof}
\begin{rem}
Equation \eqref{eq:Stark2} is similar to other representations found
in \citep{Srivastava2003}.
\end{rem}
The ingredients from the previous results can be processed together
to give a representation of \eqref{eq:Euler-Maclaurin} which does
not involve integrals: It is the compelling advantage of the following
representation over \eqref{eq:Euler-Maclaurin} and \eqref{eq:Exp-Integral}
to replace the integral by a usual sum. As the error term can be evaluated
precisely, this makes $\zeta$ available without computing the limit
and letting $n\rightarrow\infty$ in the respective sum.
\begin{thm}
For $n\in\left\{ 1,\,3,\dots\right\} $, $p\in\left\{ 2,\,3,\dots\right\} $
and $\Delta\in\left\{ 0,1,2,\dots\right\} $ the representation 
\begin{alignat}{1}
\zeta\left(s\right)= & \sum_{j=1}^{n}\frac{1}{j^{s}}+\sum_{k=0}^{\Delta-1}\frac{{s+k-2 \choose k}}{s-1}\frac{B_{k}}{n^{s+k-1}}+\nonumber \\
 & +\frac{1}{\left(s-1\right)n^{s}}\sum_{i=\Delta}\frac{{i+s-1 \choose i+1}}{\left(n\, p\right)^{i}}\sum_{j=\Delta}^{i}{i+1 \choose j}B_{j}\frac{\left(n\left(p-1\right)\right)^{i+1-j}}{p^{s}-p^{1-j}}\label{eq:Darstellung}
\end{alignat}
holds true, which represents a converging series at least for $\Re s>1-\Delta$.\end{thm}
\begin{rem}
Notice that \eqref{eq:Darstellung} even holds for $n=1$.
\end{rem}

\begin{rem}
The order of summation \emph{must} \emph{not} be changed in \eqref{eq:Darstellung}:
As formally 
\begin{alignat*}{1}
\sum_{i=\Delta}\frac{{i+s-1 \choose i}}{\left(i+1\right)\left(n\, p\right)^{i}}\sum_{j=\Delta}^{i} & {i+1 \choose j}B_{j}\frac{\left(n\left(p-1\right)\right)^{i+1-j}}{p^{s}-p^{1-j}}=\\
 & =\sum_{j=\Delta}\frac{B_{j}}{p^{s}-p^{1-j}}\sum_{i=j}^{\infty}{i+1 \choose j}\frac{{i+s-1 \choose i}}{\left(i+1\right)\left(n\, p\right)^{i}}\left(n\left(p-1\right)\right)^{i+1-j}\\
 & =\sum_{j=\Delta}\frac{B_{j}}{n^{j-1}}\frac{{s+j-2 \choose j}}{s-1},
\end{alignat*}
which characterizes the asymptotic behavior, but this series does
\emph{not }converge.
\end{rem}

\begin{rem}
For $p\rightarrow1$ and (formally) interchanging this limit with
the sum the (non-converging) identity 
\begin{align*}
\zeta\left(s\right)= & \sum_{j=1}^{n-1}\frac{1}{j^{s}}+\sum_{k=0}^{\infty}\frac{{1-s \choose k}}{s-1}\frac{B_{k}}{n^{s+k-1}}
\end{align*}
is obtained; and for $p\rightarrow\infty$, the same argument gives
\begin{align*}
\zeta\left(s\right)= & \sum_{j=1}^{n-1}\frac{1}{j^{s}}+\sum_{k=0}^{\Delta}\frac{{1-s \choose k}}{s-1}\frac{B_{k}}{n^{s+k-1}}.
\end{align*}
\end{rem}
\begin{proof}
As for the proof notice that 
\begin{alignat*}{1}
\sum_{k=0}^{\Delta-1} & \frac{{1-s \choose k}}{s-1}\frac{B_{k}}{n^{s+k-1}}+\frac{1}{n^{s}}\sum_{i=\Delta}\frac{{i+s-1 \choose i}}{\left(i+1\right)\left(n\, p\right)^{i}}\sum_{j=\Delta}^{i}{i+1 \choose j}B_{j}\frac{\left(n\left(p-1\right)\right)^{i+1-j}}{p^{s}-p^{1-j}}\\
= & \sum_{k=0}^{\Delta-1}\frac{{1-s \choose k}}{s-1}\frac{B_{k}}{n^{s+k-1}}-\frac{1}{n^{s}}\sum_{j=0}^{\Delta-1}\frac{B_{j}}{\left(n\left(p-1\right)\right)^{j-1}}\frac{1}{p^{s}-p^{1-j}}\sum_{i=j}^{\infty}\frac{{i+s-1 \choose i}}{i+1}{i+1 \choose j}\left(1-\frac{1}{p}\right)^{i}\\
 & +\frac{1}{n^{s}}\sum_{i=0}\frac{{i+s-1 \choose i}}{\left(i+1\right)\left(n\, p\right)^{i}}\sum_{j=\Delta}^{i}{i+1 \choose j}B_{j}\frac{\left(n\left(p-1\right)\right)^{i+1-j}}{p^{s}-p^{1-j}}\\
= & \left(*\right).
\end{alignat*}
Next we have for $p>1$ the geometric series 
\begin{alignat*}{1}
\frac{1}{n^{s}\left(n\, p\right)^{i}}\cdot & \frac{\left(n\left(p-1\right)\right)^{i+1-j}}{p^{s}-p^{1-j}}=\frac{\left(p-1\right)^{i+1-j}}{n^{s+j-1}p^{s+i}}\frac{1}{1-\frac{1}{p^{s-1+j}}}\\
 & =\frac{\left(p-1\right)^{i+1-j}}{n^{s+j-1}p^{s+i}}\sum_{k=0}\frac{1}{\left(p^{s-1+j}\right)^{k}}==\sum_{k=0}\frac{\left(n\left(p-1\right)p^{k}\right)^{i+1-j}}{\left(n\, p^{k+1}\right)^{s+i}},
\end{alignat*}
and thus 
\begin{alignat*}{1}
\left(*\right)= & \sum_{k=0}^{\Delta-1}\frac{{1-s \choose k}}{s-1}\frac{B_{k}}{n^{s+k-1}}-\sum_{j=0}^{\Delta-1}\frac{B_{j}}{n^{s+j-1}}\frac{1}{\left(p-1\right)^{j-1}}\frac{1}{p^{s}-p^{1-j}}\sum_{i=j}^{\infty}\frac{{i+s-1 \choose i}}{i+1}{i+1 \choose j}\left(1-\frac{1}{p}\right)^{i}\\
 & +\sum_{i=0}\frac{{i+s-1 \choose i}}{i+1}\sum_{k=0}\frac{1}{\left(n\, p^{k+1}\right)^{s+i}}\sum_{j=0}^{i}{i+1 \choose j}B_{j}\left(n\,\left(p-1\right)p^{k}\right)^{i+1-j}.
\end{alignat*}
The inner series has an explicit evaluation, as 
\[
\sum_{i=j}^{\infty}\frac{{i+s-1 \choose i}}{i+1}{i+1 \choose j}\left(1-\frac{1}{p}\right)^{i}=\left(-1\right)^{j}\left(p-1\right)^{j-1}\left(p^{s}-p^{1-j}\right)\frac{{1-s \choose j}}{s-1}.
\]
For this and the fact that $B_{2j+1}=0$ (except for $B_{1}=-\frac{1}{2}$)
the first two sums collapse to $\frac{1}{n^{s}}$, and thus 
\begin{align*}
\left(*\right)= & \frac{1}{n^{s}}+\sum_{i=0}{i+s-1 \choose i}\sum_{k=0}\frac{1}{\left(n\, p^{k+1}\right)^{s+i}}\cdot\frac{1}{i+1}\sum_{j=0}^{i}{i+1 \choose j}B_{j}\left(n\left(p-1\right)p^{k}\right)^{i+1-j}.
\end{align*}
Next observe that the inner sum is just 
\[
\frac{1}{i+1}\sum_{j=0}^{i}{i+1 \choose j}B_{j}\left(n\left(p-1\right)p^{k}\right)^{i+1-j}=\sum_{j=0}^{n\left(p-1\right)p^{k}-1}j^{i}
\]
 (cf. \eqref{eq:Faulhaber}), thus 

\begin{alignat*}{1}
\left(*\right)= & \frac{1}{n^{s}}+\sum_{i=0}{i+s-1 \choose i}\sum_{k=0}\frac{1}{\left(n\, p^{k+1}\right)^{s+i}}\cdot\sum_{j=0}^{n\left(p-1\right)p^{k}-1}j^{i}\\
 & =\frac{1}{n^{s}}+\sum_{i=0}{i+s-1 \choose i}\sum_{k=0}\frac{1}{\left(n\, p^{k+1}\right)^{s}}\cdot\sum_{j=0}^{n\left(p-1\right)p^{k}-1}\left(\frac{j}{n\, p^{k+1}}\right)^{i}.
\end{alignat*}
In this situation we may evaluate the power series and rearrange the
resulting terms as 
\begin{alignat*}{1}
\left(*\right)= & \frac{1}{n^{s}}+\sum_{k=0}\frac{1}{\left(n\, p^{k+1}\right)^{s}}\cdot\sum_{j=0}^{n\left(p-1\right)p^{k}-1}\left(\frac{1}{1-\frac{j}{n\, p^{k+1}}}\right)^{s}\\
 & =\frac{1}{n^{s}}+\sum_{k=0}\sum_{j=0}^{n\left(p-1\right)p^{k}-1}\left(\frac{1}{n\, p^{k+1}-j}\right)^{s}\\
 & =\frac{1}{n^{s}}+\sum_{k=0}\sum_{j=1+n\, p^{k}}^{n\, p^{k+1}}\left(\frac{1}{j}\right)^{s}=\sum_{j=n}^{\infty}\frac{1}{j^{s}}\\
 & =\zeta\left(s\right)-\sum_{j=1}^{n-1}\frac{1}{j^{s}},
\end{alignat*}
which finally completes the proof.

To analyze convergence it should be mentioned that 
\begin{align*}
\frac{1}{n^{s}}\sum_{i=\Delta} & \frac{{s-1+i \choose i}}{\left(i+1\right)\left(n\, p\right)^{i}}\sum_{j=\Delta}^{i}{i+1 \choose j}B_{j}\frac{\left(n\left(p-1\right)\right)^{i+1-j}}{p^{s}-p^{1-j}}\\
= & \frac{1}{n^{s}}\sum_{i=\Delta}{s-1+i \choose i}*\\
 & \sum_{k=0}\frac{1}{\left(n\, p^{k+1}\right)^{s+i}}\left\{ \sum_{j=0}^{n\left(p-1\right)p^{k}-1}j^{i}-\frac{1}{1+i}\sum_{j=0}^{\Delta-1}{i+1 \choose j}\cdot B_{j}\cdot\left(n\left(p-1\right)p^{k}\right)^{i+1-j}\right\} .
\end{align*}
Recall from \eqref{eq:Euler-Maclaurin} and \eqref{eq:Faulhaber}
that $\sum_{j=1}^{n-1}j^{i}=\mathcal{O}\left(\frac{n^{i+1}}{i+1}\right)$,
and further that 
\[
\sum_{j=1}^{n-1}j^{i}-\frac{1}{1+i}\sum_{j=0}^{\Delta-1}{i+1 \choose j}\cdot B_{j}\cdot n^{i+1-j}=\mathcal{O}\left(n^{i+1-\Delta}\right).
\]
Whence convergence can be described by 
\begin{alignat*}{1}
\frac{1}{\left(n\, p^{k+1}\right)^{s+i}} & \left\{ \sum_{j=0}^{n\left(p-1\right)p^{k}-1}j^{i}-\frac{1}{1+i}\sum_{j=0}^{\Delta-1}{i+1 \choose j}\cdot B_{j}\cdot\left(n\left(p-1\right)p^{k}\right)^{i+1-j}\right\} \\
= & \mathcal{O}\left(\frac{\left(n\left(p-1\right)p^{k}\right)^{i+1-\Delta}}{\left(n\, p^{k+1}\right)^{s+i}}\right)\\
= & \mathcal{O}\left(n^{1-\Delta-s}\frac{\left(p-1\right)^{i+1-\Delta}}{p^{s+i}}p^{k\left(1-\Delta-s\right)}\right)\\
= & \mathcal{O}\left(\left(1-\frac{1}{p}\right)^{i}\frac{1}{p^{k\left(s-1+\Delta\right)}}\right),
\end{alignat*}
which shows that the series converge at least for $\Re s>1-\Delta$.
\end{proof}
Woon's formula (cf. \citep{Borwein2000,Woon1998}) is to some (minor)
extend close to the representation above, it involves a double sum
as our representation \eqref{eq:Darstellung}.

\section{Parametrization by Complete Elliptic Integrals}

Riemann's Zeta function has the representation 
\begin{equation}
\zeta\left(s\right)\frac{\Gamma\left(\tfrac{s}{2}\right)}{\pi^{\frac{s}{2}}}=\int_{0}^{\infty}t^{\frac{s}{2}-1}\frac{1}{2}\left(\vartheta_{3}\left(0;it\right)-1\right)\mathrm{d}t,\label{eq:theta}
\end{equation}
$\theta\left(t\right)=\vartheta_{3}\left(0;it\right)=\sum_{n\in\mathbb{Z}}q^{n^{2}}=1+2\sum_{n=1}q^{n^{2}}$
being Jacobi's elliptic theta function (we choose $q=e^{-i\pi t}$
here for convenience; \citep{Edwards1974} uses $\Psi$ for $\Psi\left(t\right)=\frac{1}{2}\left(\vartheta_{3}\left(0;it\right)-1\right)=\sum_{n=1}e^{-n^{2}\pi x}$).
The latter representation \eqref{eq:theta} deduces immediately from
\[
\Gamma\left(\frac{s}{2}\right)=n^{s}\pi^{\frac{s}{2}}\cdot\int_{0}^{\infty}t^{\frac{s}{2}-1}e^{-n^{2}\pi t}\mathrm{d}t
\]
and holds true for $\Re s>1$. In order to get a formula for \emph{all}
complex numbers $s$ one may apply integration by parts twice for
\begin{equation}
s\left(s-1\right)\zeta\left(s\right)\frac{\Gamma\left(\frac{s}{2}\right)}{\pi^{\frac{s}{2}}}=\int_{0}^{\infty}t^{\frac{s}{2}-1}\tilde{\theta}\left(t\right)\mathrm{d}t,\label{eq:thetaw}
\end{equation}
where 
\[
\tilde{\theta}\left(t\right)=2t^{2}\,\theta^{\prime\prime}\left(t\right)+3t\,\theta^{\prime}\left(t\right)=\sum_{n=1}\left(4\left(n^{2}\pi t\right)^{2}-6n^{2}\pi t\right)e^{-n^{2}\pi t}.
\]
It can be proved directly that $\tilde{\theta}$ inherits the property
\begin{equation}
\tilde{\theta}\left(t\right)=\frac{1}{\sqrt{t}}\tilde{\theta}\left(\frac{1}{t}\right)\label{eq:einsdurch}
\end{equation}
 from $\theta$, but additionally 
\[
0=\lim_{x\rightarrow\infty}\tilde{\theta}^{\left(k\right)}\left(x\right)=\lim_{x\rightarrow0}\tilde{\theta}^{\left(k\right)}\left(x\right)
\]
for all derivatives including the function $\tilde{\theta}$ itself
($k=0,1\dots$) holds true, such that the next integral \eqref{eq:thetaw}
defines an analytic function for all values of $s\in\mathbb{C}$. 

The symmetric version (for $e^{\frac{x}{4}}\tilde{\theta}\left(e^{x}\right)=e^{-\frac{x}{4}}\tilde{\theta}\left(e^{-x}\right)$,
again by \eqref{eq:einsdurch}) of \eqref{eq:thetaw} reads
\begin{equation}
s\left(s-1\right)\zeta\left(s\right)\frac{\Gamma\left(\frac{s}{2}\right)}{\pi^{\frac{s}{2}}}=\int_{-\infty}^{\infty}e^{\frac{x}{2}\left(s-\frac{1}{2}\right)}\cdot e^{\frac{x}{4}}\tilde{\theta}\left(e^{x}\right)\mathrm{d}x,\label{eq:thetaw-2}
\end{equation}
and it is an interesting observation that this is the continuous Fourier
transform of $x\mapsto e^{\frac{x}{4}}\tilde{\theta}\left(e^{x}\right)$:
\[
s\left(s-1\right)\zeta\left(s\right)\frac{\Gamma\left(\frac{s}{2}\right)}{\pi^{\frac{s}{2}}}=\sqrt{2\pi}\cdot\mathcal{F}\left(e^{\frac{x}{4}}\tilde{\theta}\left(e^{x}\right)\right)\left(\tfrac{i}{2}\left(s-\tfrac{1}{2}\right)\right),
\]
where $\mathcal{F}\left(f\right)\left(k\right)=\frac{1}{\sqrt{2\pi}}\int_{-\infty}^{\infty}e^{ikx}f(x)\mathrm{d}x$
is the Fourier transform of $f$.

The property \eqref{eq:einsdurch} may be used to easily deduce Riemann's
functional equation. The function decays rather quickly to $0$, for
$e^{\frac{x}{4}}\tilde{\theta}\left(e^{x}\right)\sim4\pi^{2}e^{\frac{9}{4}x}e^{-\pi e^{-x}}$
as $x$ tends to infinity. 

A good proxy of \eqref{eq:thetaw} in terms of a polynomial in $s$
with convergence as $n\rightarrow\infty$ is given by the asymptotic
properties of the Hermite polynomial $H_{n}$ (cf. \citep{Polya1927})
and Laguerre polynomial $L_{n}$, 
\begin{align}
\frac{\left(-1\right)^{n}\cdot n!}{(2n)!} & \int_{-\infty}^{\infty}H_{2n}\left(x\cdot\frac{s-\frac{1}{2}}{4i\sqrt{n}}\right)\cdot e^{\frac{x}{4}}\tilde{\theta}\left(e^{x}\right)\mathrm{d}x=\label{eq:Hermite}\\
= & \frac{1}{{n-\tfrac{1}{2} \choose n}}\int_{-\infty}^{\infty}L_{n}^{\left(-\frac{1}{2}\right)}\left(\frac{x^{2}}{n}\cdot\left(\frac{s-\frac{1}{2}}{4i}\right)^{2}\right)\cdot e^{\frac{x}{4}}\tilde{\theta}\left(e^{x}\right)\mathrm{d}x\nonumber \\
= & \frac{1}{{n-\tfrac{1}{2} \choose n}}\int_{-\infty}^{\infty}L_{n}^{\left(-\frac{1}{2}\right)}\left(\frac{x^{2}}{16n}\cdot\left(\frac{s-\frac{1}{2}}{i}\right)^{2}\right)\cdot e^{\frac{x}{4}}\tilde{\theta}\left(e^{x}\right)\mathrm{d}x\nonumber \\
 & \xrightarrow[n\to\infty]{}s\left(s-1\right)\zeta\left(s\right)\frac{\Gamma\left(\frac{s}{2}\right)}{\pi^{\frac{s}{2}}}.\nonumber 
\end{align}

The symmetric integrand is maintained by the following form in terms
of complete Elliptic Integrals of the first and second kind.
\begin{thm}
For any $s\in\mathbb{C}$ 
\begin{equation}
s\left(s-1\right)\zeta\left(s\right)\frac{\Gamma\left(\frac{s}{2}\right)}{\pi^{\frac{s}{2}}}=\int_{0}^{1}\left(\frac{\mathbf{K}\left(1-m\right)}{\mathbf{K}\left(m\right)}\right)^{\frac{1}{2}\left(s-\frac{1}{2}\right)}U\left(m\right)\mathrm{d}m,\label{eq:thetaw-1}
\end{equation}
where $U$ is the symmetric function 
\[
U\left(m\right)=\frac{\left(\mathbf{K}\left(m\right)\mathbf{K}\left(1-m\right)\right)^{\frac{1}{4}}}{\pi\sqrt{2\pi}\, m\left(1-m\right)}\left(\begin{array}{r}
3\left(1-m\right)\cdot\mathbf{K}\left(m\right)\mathbf{E}\left(1-m\right)\\
+3m\cdot\mathbf{E}\left(m\right)\mathbf{K}\left(1-m\right)\\
-3\cdot\mathbf{E}\left(m\right)\mathbf{E}\left(1-m\right)\\
-m\left(1-m\right)\cdot\mathbf{K}\left(m\right)\mathbf{K}\left(1-m\right)
\end{array}\right).
\]
Here, $\mathbf{K}\left(m\right)=\int_{0}^{\frac{\pi}{2}}\frac{\mathrm{d}x}{\sqrt{1-m\sin^{2}x}}=\frac{\pi}{2}\,_{2}F_{1}\left(\frac{1}{2},\frac{1}{2},1,m\right)$
is the \emph{complete elliptic integral of the first kind}, and $\mathbf{E}\left(m\right)=\int_{0}^{\frac{\pi}{2}}\sqrt{1-m\sin^{2}x}\mathrm{d}x=\frac{\pi}{2}\,_{2}F_{1}\left(\frac{1}{2},-\frac{1}{2},1,m\right)$
the \emph{complete elliptic integral of the second kind}%
\footnote{\emph{Elliptic integrals are often defined with the parameter $m=k^{2}$
instead; in the present context, for symmetry, it is more convenient
to use the parameter $m$. }%
}, which describes the circumference of the ellipse with eccentricity
$\sqrt{m}$.

\begin{figure}
\includegraphics[width=0.33\textwidth]{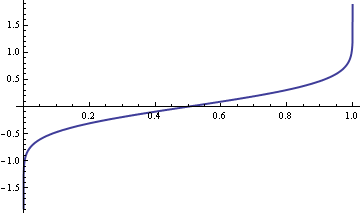} \includegraphics[width=0.33\textwidth]{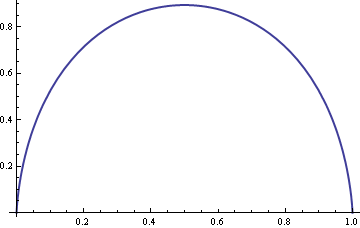}
\includegraphics[width=0.33\textwidth]{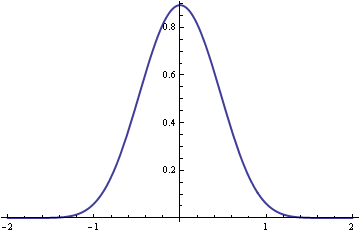}

\caption{Plot of the symmetric functions $m\mapsto x\left(m\right)$, $m\mapsto\tilde{U}\left(m\right)$
and $x\mapsto e^{\frac{x}{4}}\tilde{\vartheta}\left(e^{x}\right)$.}

\end{figure}

\end{thm}
We shall use the abbreviations $\mathbf{K}=\mathbf{K}\left(m\right)$,
$\mathbf{K}^{\prime}=\mathbf{K}\left(1-m\right)$, $\mathbf{E}=\mathbf{E}\left(m\right)$
and $\mathbf{E}^{\prime}=\mathbf{E}\left(1-m\right)$ in the sequel.

A power series analysis of $U$ gives that 
\[
U\left(m\right)=\frac{1}{4\,\sqrt[4]{\pi}}\left(\ln\frac{16}{m}-\frac{3}{2}\right)\left(\ln\frac{16}{m}\right)^{\frac{1}{4}}+o\left(1\right)
\]
and 
\[
\frac{\mathbf{K}\left(1-m\right)}{\mathbf{K}\left(m\right)}=\frac{1}{\pi}\ln\frac{16}{m}+o\left(1\right).
\]
As 
\[
\int_{0}^{x}\left(\ln\frac{16}{m}\right)^{\alpha}\mathrm{d}m=16\cdot\Gamma\left(\alpha+1,\ln\frac{16}{x}\right)\approx x\left(\ln\frac{16}{x}\right)^{\alpha}\left(1+\frac{\alpha}{\ln\frac{16}{x}}+\dots\right)
\]
is integrable for any $\alpha$ with a result expressed by the upper
incomplete Gamma function $\Gamma$. It is thus evident that the representation
\eqref{eq:thetaw-1} converges for all values of $s$.
\begin{proof}
For the nome $q\left(m\right):=e^{-\pi\frac{\mathbf{K}\left(1-m\right)}{\mathbf{K}\left(m\right)}}$
the q-series identity 
\begin{equation}
1+2\sum_{n=1}q\left(m\right)^{n^{2}}=\sqrt{\frac{2}{\pi}\mathbf{K}\left(m\right)}\label{eq:qK}
\end{equation}
holds true. Observe that 
\[
\frac{\mathrm{d}}{\mathrm{d}m}\log q\left(m\right)=-\pi\frac{\mathrm{d}}{\mathrm{d}m}\frac{\mathbf{K}\left(1-m\right)}{\mathbf{K}\left(m\right)}=-\pi\frac{\mathbf{E}\mathbf{K}^{\prime}+\mathbf{K}\mathbf{E}^{\prime}-\mathbf{K}\mathbf{K}^{\prime}}{2m\left(1-m\right)\mathbf{K}^{2}}=\frac{\pi^{2}}{4m\left(1-m\right)\mathbf{K}^{2}}
\]
 due to the well-known derivatives for $\frac{\mathrm{d}}{\mathrm{d}m}\mathbf{K}=\frac{\mathbf{E}-\left(1-m\right)\mathbf{K}}{2m\left(1-m\right)}$
and $\frac{\mathrm{d}}{\mathrm{d}m}\mathbf{E}=\frac{\mathbf{E}-\mathbf{K}}{2m}$
and Legendre's relation (cf. \citep{Erdelyi1953} or \citep{Gradsteyn2007}). 

Taking the derivative of the formula \eqref{eq:qK} reveals the closed
form 
\[
\sum_{n=1}q\left(m\right)^{n^{2}}n^{2}\pi=\sqrt{\frac{2}{\pi}\mathbf{K}\left(m\right)}\frac{\mathbf{K}}{2\pi}\left(\mathbf{E}-\left(1-m\right)\mathbf{K}\right)
\]
after simplification, and proceeding in the same way further 
\[
\sum_{n=1}q\left(m\right)^{n^{2}}\left(n^{2}\pi\right)^{2}=\sqrt{\frac{2}{\pi}\mathbf{K}\left(m\right)}\frac{\mathbf{K}^{2}}{2\pi^{2}}\left(3\mathbf{E}^{2}-6\left(1-m\right)\mathbf{E}\mathbf{K}+\left(3-m\right)\left(1-m\right)\mathbf{K}^{2}\right).
\]

Combining the latter equalities and after further simplifications
involving Legendre's relation again one obtains 
\[
t\left(m\right)^{-\frac{3}{4}}\cdot\tilde{\theta}\left(t\left(m\right)\right)\cdot t^{\prime}\left(m\right)=U\left(m\right),
\]
where $t\left(m\right):=-\frac{1}{\pi}\log q\left(m\right)=\frac{\mathbf{K}\left(1-m\right)}{\mathbf{K}\left(m\right)}$.
Whence 
\begin{alignat*}{1}
s\left(s-1\right)\zeta\left(s\right) & \frac{\Gamma\left(\frac{s}{2}\right)}{\pi^{\frac{s}{2}}}=\int_{0}^{\infty}t^{\frac{s}{2}-1}\tilde{\theta}\left(t\right)\mathrm{d}t\\
 & =\int_{0}^{1}t\left(m\right)^{\frac{s}{2}-1}\tilde{\theta}\left(t\left(m\right)\right)t^{\prime}\left(m\right)\mathrm{d}m=\int_{0}^{1}t\left(m\right)^{\frac{s-\frac{1}{2}}{2}}t\left(m\right)^{-\frac{3}{4}}\tilde{\theta}\left(t\left(m\right)\right)t^{\prime}\left(m\right)\mathrm{d}m\\
 & =\int_{0}^{1}\left(\frac{\mathbf{K}\left(1-m\right)}{\mathbf{K}\left(m\right)}\right)^{\frac{1}{2}\left(s-\frac{1}{2}\right)}U\left(m\right)\mathrm{d}m,
\end{alignat*}
which is the desired assertion. 
\end{proof}
The results found can be combined in the following way. 

Recall the inverse elliptic nome function $q^{-1}\left(q\right):=\frac{\vartheta_{2}\left(q\right)^{4}}{\vartheta_{3}\left(q\right)^{4}}$
(cf. \citep[Chapter XXI, p. 486]{Whittaker1902})%
\footnote{This is called the \emph{problem of inversion}, cf. \citep[volume II, p. 362]{Erdelyi1953}.%
}; with this choice $x^{-1}\left(x\right)=q^{-1}\left(e^{-\pi e^{-x}}\right)$
is the inverse function of $x\left(m\right):=\log\frac{\mathbf{K}\left(m\right)}{\mathbf{K}\left(1-m\right)}$. 

This can be used to rewrite \eqref{eq:thetaw-2} in different variants
as 
\begin{alignat*}{1}
s\left(s-1\right)\zeta\left(s\right)\frac{\Gamma\left(\frac{s}{2}\right)}{\pi^{\frac{s}{2}}} & =\int_{-\infty}^{\infty}e^{\frac{1}{2}\left(s-\frac{1}{2}\right)x}\frac{U\left(x^{-1}\left(x\right)\right)}{x^{\prime}\left(x^{-1}\left(x\right)\right)}\mathrm{d}x\\
 & =\int_{-\infty}^{\infty}e^{\frac{1}{2}\left(s-\frac{1}{2}\right)x}\cdot\tilde{U}\left(\left(\frac{\vartheta_{2}\left(e^{-\pi e^{-x}}\right)}{\vartheta_{3}\left(e^{-\pi e^{-x}}\right)}\right)^{4}\right)\mathrm{d}x\\
 & =\int_{0}^{\infty}t^{\frac{s-\frac{1}{2}-2}{2}}\cdot\tilde{U}\left(\frac{\vartheta_{2}^{4}\left(e^{-\pi t}\right)}{\vartheta_{3}^{4}\left(e^{-\pi t}\right)}\right)\mathrm{d}t\\
 & =\int_{0}^{1}\left(-\frac{\log q}{\pi}\right)^{\frac{s-\frac{1}{2}-2}{2}}\frac{1}{q\pi}\tilde{U}\left(\frac{\vartheta_{2}^{4}\left(q\right)}{\vartheta_{3}^{4}\left(q\right)}\right)\mathrm{d}q
\end{alignat*}
where
\[
\tilde{U}\left(m\right)=4\frac{\left(\mathbf{K}\left(m\right)\mathbf{K}\left(1-m\right)\right)^{\frac{5}{4}}}{\pi^{2}\sqrt{2\pi}}\left(\begin{array}{r}
3\left(1-m\right)\cdot\mathbf{K}\left(m\right)\mathbf{E}\left(1-m\right)\\
+3m\cdot\mathbf{E}\left(m\right)\mathbf{K}\left(1-m\right)\\
-3\cdot\mathbf{E}\left(m\right)\mathbf{E}\left(1-m\right)\\
-m\left(1-m\right)\cdot\mathbf{K}\left(m\right)\mathbf{K}\left(1-m\right)
\end{array}\right).
\]
Introducing $\mathbf{B}\left(m\right):=\frac{\mathbf{E}-\left(1-m\right)\mathbf{K}}{m}=\int_{0}^{\frac{\pi}{2}}\frac{\cos^{2}\phi}{\sqrt{1-m\sin^{2}\phi}}\mathrm{d}\phi=\frac{\pi}{4}\,_{2}F_{1}\left(\frac{1}{2},\frac{1}{2},2,m\right)$
this function rewrites as 
\[
\tilde{U}\left(m\right)=\sqrt{8\pi}\frac{m\left(1-m\right)}{\pi^{3}}\left(\mathbf{K}^{\prime}\mathbf{K}\right)^{\frac{5}{4}}\left(2\mathbf{K}^{\prime}\mathbf{K}-3\mathbf{B}^{\prime}\mathbf{B}\right).
\]

Comparing the identity with \eqref{eq:thetaw-2} we finally obtain
that $t^{\frac{1}{4}}\tilde{\vartheta}\left(t\right)=\tilde{U}\left(\left(\frac{\vartheta_{2}\left(e^{-\pi t}\right)}{\vartheta_{3}\left(e^{-\pi t}\right)}\right)^{4}\right)=\tilde{U}\left(\left(\frac{\vartheta_{4}\left(e^{-\pi t}\right)}{\vartheta_{3}\left(e^{-\pi t}\right)}\right)^{4}\right)$,
or 
\[
e^{\frac{x}{4}}\tilde{\theta}\left(e^{x}\right)=\tilde{U}\left(\left(\frac{\vartheta_{2}\left(e^{-\pi e^{-x}}\right)}{\vartheta_{3}\left(e^{-\pi e^{-x}}\right)}\right)^{4}\right)=\tilde{U}\left(q^{-1}\left(e^{-\pi e^{-x}}\right)\right)
\]
by symmetry and \eqref{eq:einsdurch}. The inversion moreover has
the well-known (cf. \citep[Volume II, p. 362]{Erdelyi1953}, correcting
the misprints there) expansion $\frac{\vartheta_{4}\left(q\right)}{\vartheta_{3}\left(q\right)}=\left(\prod_{n=0}\frac{1-q^{2n+1}}{1+q^{2n+1}}\right)^{2}$
as an infinite product.

The power series expansion 
\begin{alignat*}{1}
s\left(s-1\right)\zeta\left(s\right)\frac{\Gamma\left(\frac{s}{2}\right)}{\pi^{\frac{s}{2}}} & =\sum_{k=0}\frac{\left(s-\frac{1}{2}\right)^{2k}}{(2k)!4^{k}}\int_{-\infty}^{\infty}x^{2k}\cdot e^{\frac{x}{4}}\tilde{\theta}\left(e^{x}\right)\mathrm{d}x\\
 & =\sum_{k=0}\frac{\left(s-\frac{1}{2}\right)^{2k}}{(2k)!4^{k}}\int_{-\infty}^{\infty}x^{2k}\cdot\tilde{U}\left(\left(\frac{\vartheta_{2}\left(e^{-\pi e^{-x}}\right)}{\vartheta_{3}\left(e^{-\pi e^{-x}}\right)}\right)^{4}\right)\mathrm{d}x
\end{alignat*}
for even and positive coefficients then is just a by-product. The
convergence of \eqref{eq:Hermite}, however, is faster.

\section{Summary And Acknowledgment}

\subsection{Summary}

Various algorithms to evaluate the Riemann $\zeta$ function start
with the Euler-Maclaurin summation formula, as it properly describes
the asymptotics of the series representation. The closed forms available
for the error bound typically involve integrals. The method proposed
here generalizes other methods and gives the precise error bound as
sum.

\subsection{Acknowledgment}

We would like to thank the reviewers for their significant contribution
to improve the paper.

The plots have been produced by use of Mathematica.
\begin{figure}[h]
\centering{}\includegraphics[scale=0.3]{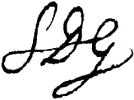}
\end{figure}

\section{References}

\bibliographystyle{alpha}
\bibliography{../../Literatur/LiteraturAlois}

\end{document}